\documentclass{article}

\usepackage[utf8]{inputenc}
\usepackage[russian,english]{babel}
\usepackage{amsthm,amssymb,amsmath,amsfonts,mathrsfs}
\usepackage{graphicx}

\newtheorem*{theorem}{Theorem}

\usepackage{titling}
\title{Square-tiled surfaces and curves over number fields}

\author{George B.~Shabat
\thanks{RSUH (Russian State University for the Humanities), Moscow, Russia.
\newline e-mail: {\tt george.shabat@gmail.com}}}

\date{}

\providecommand{\keywords}[1]{{\bf Keywords:} #1}

\begin{document}

\maketitle


\begin{abstract}

The paper presents an analog of the old result by the author and V. Voevodsky, according to which a Riemann surface admits a conformal structure, defined by an equilateral triangulation, if and only if the corresponding algebraic curve can be defined over the field of the algebraic numbers; the similar result is obtained for the square-tiled surfaces.
\end{abstract}

\section{Introduction}

This paper presents an analog of the old result \cite {ShVo1989}, according to which a Riemann surface admits a conformal structure, defined by an equilateral triangulation, if and only if the corresponding algebraic curve can be defined over the field of the algebraic numbers. Here we establish the similar result, replacing the equilateral triangles by squares; one of the important differences between this and the previous paper is that in 1980-s the piece-wise euclidean surfaces were somewhat exotic while nowadays the {\it square-tiled} ones or {\it origamis} are studied intensively during the last decades; see, e.g., \cite{Zor2006}.

Among the immediate applications we obtain the possibility of explicitly writing down the defining equations of some curves traditionally defined by the transcendental means. {\bf
 Hopefully, the strong statistical results concerning the square-tiled surfaces can promote the understanding of the typical behavior of dessins d'enfants with many edges -- e.g., the distribution of sizes of their Galois orbits.}

\section{Preliminaries}

We are going to work with the {\it piecewise-euclidean} compact connected orientable closed surfaces. There are two equivalent ways to define this structure. One of them is purely metrical and is defined by providing a surface with a metric that is riemannian away from a finite set and is {\it flat} away from this set; briefly the surfaces themselves are also called {\it flat}, see \cite{Zor2006}.

We take another approach, considering surfaces with some combinatorial structure, a {\it dessin d'enfant} (see \cite{Gro1984}), and then supply a dessin with an additional metric structure. Topologically a dessin is a triple
$$
\mathcal{D}=(\mathbf{X}_2\supset\mathbf{X}_1\supset\mathbf{X}_0),
$$
where $\mathbf{X}_2$ is an oriented compact connected surface, $\mathbf{X}_0$ is a finite set and  $\mathbf{X}_i\setminus \mathbf{X}_{i-1}$ is homeomorphic to a finite disjpoint union of $i$-cells for $i\in\{1,2\}$.

Informally a surface is just subdivided into the union of metric polygons. A convenient {\it cartographic} formalism was suggested in \cite{Gro1984} (though well-known long before -- %to combinatorial people, 
not to Grothendieck...). The combinatorics of a dessin $\mathcal{D}$ is fixed by the set of the {\it directed edges} $\vec{E}(\mathcal{D})$, acted upon by the {\it oriented cartographic group}
$$
\mathcal{C}_2^+:=\langle \rho_0, \rho_1,\rho_2\mid
\rho_1^2=\rho_2 \rho_1\rho_0=1\rangle;
$$
$\rho_0$ rotates the directed edge contrary-clockwise, $\rho_1$ reverses the directed edge and $\rho_2$ moves it along the boundary of the cell that lies \underline{to the left} of it. The sets of vertexes, edges and faces of $\mathcal{D}$ are restored by
$$
V(\mathcal{D}):=\frac{\vec{E}(\mathcal{D})}{\langle\rho_0\rangle},
$$
%%%%%%%%
$$
E(\mathcal{D}):=\frac{\vec{E}(\mathcal{D})}{\langle\rho_1\rangle}
$$
and
$$
F(\mathcal{D}):=\frac{\vec{E}(\mathcal{D})}{\langle\rho_2\rangle},
$$
with the incidence relations defined in terms of the subgroups of $\mathcal{C}_2^+$.

For a directed edge $\varepsilon\in\vec{E}(\mathcal{D})$ we shall use the Kontsevich \cite{{Konts1992}} notations 
$[\varepsilon]_j\in\frac{\vec{E}(\mathcal{D})}{\langle\rho_j\rangle}$ for 
$j\in{0,1,2}$; thus $[\varepsilon]_0$ is the origin of a directed edge, 
$[\varepsilon]_1$ is the corresponding non-oriented edge, and $[\varepsilon]_2$ is the 2-cell ``to the left'' of it.

An additional structure that turns the surface $\mathbf{X}_2$ into the metric space is defined face-by-face; we fix it by choosing the complex coordinate $z:U\hookrightarrow\mathbb{C}$ for each %small neighborhood of a
 connected component $U$ of $\mathbf{X}_2\setminus\mathbf{X}_1$. It is assumed that the closure of each $z(U)$ is a polygon and the metric on $U$ is induced by the standard euclidean metric on $\mathbb{C}$, fixed by the \it lengths\rm
$$
\ell:\vec{E}(\mathcal{D})\longrightarrow\mathbb{R}_{>0}
$$
and \it angles\rm
$$
\phi:\vec{E}(\mathcal{D})\longrightarrow(0,2\pi)
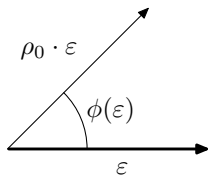
$$
\begin{center}
\includegraphics{quad_colors-101.mps}
\end{center}
The $z'$s are defined up to $z\leftarrow az+b$ with 
$a,b\in\mathbb{C}$, where $|a|=1$. We choose the canonical representatives
$z_\varepsilon$ of these coordinates, parametrized by all the $\varepsilon\in\vec{E}(\mathcal{D})$, normalized by the conditions $z_\varepsilon([\varepsilon]_0)=0$ and $z_\varepsilon([\varepsilon]_1)\subset\mathbb{R}_{\ge0}$. Then the coordinates $z_\varepsilon$ are supposed to be extended to some neighborhoods of the cells $[\varepsilon]_2$ (not containing the vertexes), and in the common domains of definition the relations
$$
z_{\rho_0\cdot\varepsilon}=\mathrm{e}^{\mathrm{i}\phi(\varepsilon)}z_\varepsilon,
$$
%%%
$$
z_{\rho_1\cdot\varepsilon}=\ell(\varepsilon)-z_\varepsilon
$$
hold.
The numbers $\ell(\varepsilon)$ and $\phi(\varepsilon)$ should satisfy some relations implying that {\it every 2-face metrically is a closed polygon}. In the simplest cases that we are going to consider these relations are obvious.

For every dessin $\mathcal{D}$ denote 
$$
\mathrm{Met}(\mathcal{D})\subset\mathbb{R}^{\vec{E}(\mathcal{D})}\times(0,2\pi)^{\vec{E}(\mathcal{D})}
$$
the set of metrical structures on the surface $\mathbf{X}_2$ on which the dessin $\mathcal{D}$ is drawn.

Then denote, as usual, by $\mathcal{DESS}$ the category of dessins and by $\mathbf{DESS}$ the set of classes of isomorphism of dessins, split this set by the genera 
$$
\mathbf{DESS}=\coprod_{g=0}^\infty\mathbf{DESS}_g.
$$
Consider for every $g\in\mathbb{N}$ the set
$$
\mathrm{PE}_g:=\bigcup_{[\mathcal{D}]\in\mathbf{DESS}_g}\mathrm{Met}(\mathcal{D})
$$
of {\it piecewise-euclidean} structures on the surface of genus $g$. It carries the structure of infinite-dimensional real manifold with boundary, but we will not use it.

\section{Main theorem}

In the previous section we have associated to every metric $\mu\in\mathrm{PE}_g$ the set of local holomorphic coordinates $\{z_\varepsilon\}$, thus defining the mapping
$$
\mathbf{compl}_g:\mathrm{PE}_g\longrightarrow\mathcal{M}_g(\mathbb{C})
$$
that transforms the metrical structures to the complex ones. It is surjective, since on every Riemann surface a \it Strebel differential \rm can be found (see \cite{Streb1984}) and therefore it can be subdivided into rectangles.

We are interested, however, in the $\mathbf{compl}_g$-images of some countable subsets of $\mathrm{PE}_g$. Introduce the set
$\Delta_g\subset\mathrm{PE}_g$ of \it equilateral triangulations \rm (defined by $\ell\equiv1$, $\phi\equiv\frac{\pi}{3}$) and
the set $\square_g\subset\mathrm{PE}_g$ of \it square tilings \rm (here $\ell\equiv1$, $\phi\equiv\frac{\pi}{2}$). It was proved in \cite{ShVo1989} that $\mathbf{compl}_g(\Delta_g)=\mathcal{M}_g(\overline{\mathbb{Q}})$. The goal of the present paper is to prove the similar result for the square tilings.

\begin{theorem}
The equality $\mathbf{compl}_g(\square_g)=\mathcal{M}_g(\overline{\mathbb{Q}})$ holds. In other words, a complex curve $\mathbf{X}$ can be defined over the field of algebraic numbers if and only if its complex structure can be defined by the piece-wise euclidean metric, in which the Riemann surface of $\mathbf{X}$ is the union of squares.
\end{theorem}

\bf Proof. \rm \underline{The ``if'' part}. Separate each square by its diagonals into four right isosceles triangles and color them like this:

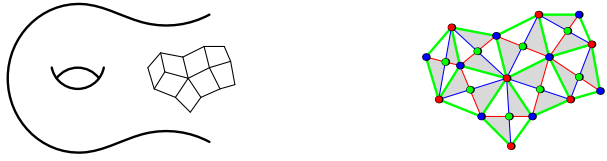
\begin{center}
\includegraphics{quad_colors-102.mps}
\end{center}
Namely, introduce on the surface $\mathbf{X}$ a \it tricolored dessin \rm (see \cite{Shab2016}), i.e. 
the dessin
$\mathbf{X}=\mathbf{X}_2\supset\mathbf{X}_1\supset\mathbf{X}_0 $ endowed with a
\it coloring mapping \rm
$$
\mathrm{col}_3:\mathbf{X}_1\longrightarrow\{blue,\,green,\,red\},
$$
constant on the edges. It is demanded that
\begin{enumerate}
\item[{\bf (0)}] any vertex is incident to edges of exactly \underline{two} colors;\\
\item[{\bf (1)}] any edge has \underline{two} vertices in its closure;\\ 
\item[{\bf (2)}] any face has \underline{three} edges in its closure, colored pairwise differently.
\end{enumerate}

Then for each square with the vertices $A,B,C,D\in\mathbf{X}_0$ denote $\varepsilon$ the directed edge ``from $A$ to $B$''; in the above notations it means that $[\varepsilon]_0=A$ and $[\rho_1\cdot\varepsilon]_0=B$). The relations 
$[\rho_2\cdot\varepsilon]_0=B$, $[\rho_2^2\cdot\varepsilon]_0=C$ and $[\rho_2^3\cdot\varepsilon]_0=D$ hold, and 
 the holomorphic coordinate $z_\varepsilon$ satisfies
$$
z_\varepsilon(A)=0,\  z_\varepsilon(B)=1,\ z_\varepsilon(C)=1+\mathrm{i},\ z_\varepsilon(D)=\mathrm{i}.
$$
Denote $O$ the intersection of the diagonals of the square $ABCD$, then $z_\varepsilon(O)=\frac{1+\mathrm{i}}{2}$.\\
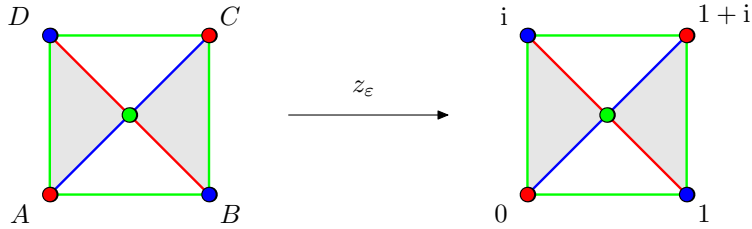
\begin{center}
\includegraphics{quad_colors-103.mps}
\end{center}
Now we can define the Belyi function $\beta:\mathbf{X}\to\mathbf{P}_1(\mathbb{C})$ demanding that its restriction to the interior of the triangle $AOD$ realizes the conformal equivalence with the lower half-plane 
$\mathcal{L}:=\{w\in\mathbb{C}\mid\mathrm{Im}(w)<0\}$
$$
\beta|_{AOD}:AOD\stackrel{\simeq}\longrightarrow\mathcal{L},
$$
satisfying the boundary conditions
$$
\beta(O)=-\infty, \beta(A)=0, \beta(D)=1.
$$
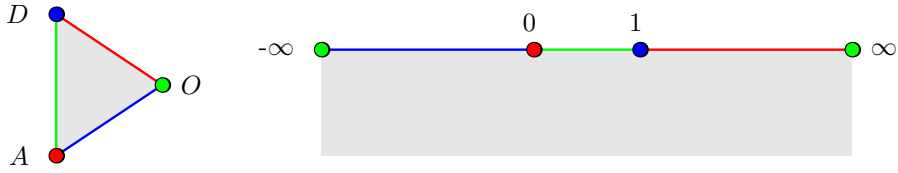
\begin{center}
\includegraphics{quad_colors-104.mps}
\end{center}
Such $\beta$ by the \it Christoffel-Schwartz \rm integral (see, e.g., \cite{Neh1982}) can be written explicitly in terms of the holomorphic
 coordinate
$z_\varepsilon$: for a point $P$ inside the considered triangle $ AOD$
$$
z(P)=\mathrm{i}\frac{\displaystyle\int\limits_0^{\beta(P)}\frac{\mathrm{d}w}{\sqrt[4]{w^3(1-w)^3}}}
{\displaystyle\int\limits_0^1\frac{\mathrm{d}w}{\sqrt[4]{w^3(1-w)^3}}}.
$$
%where $C$ is a constant chosen so that \\
It follows from the general principles of the geometric function theory that such $\beta$'s defined ``triangle-wise'' are the restrictions of the global meromorphic function on $\mathbf{X}$ that maps the white triangles to the upper half-plane and the black triangles to the lower half-plane. This function realizes the conformal map away from the vertices, so it branches only over $0,1,\infty$ and hence is a Belyi function. The assertion follows from the easy half of Belyi's theorem.

\rm \underline{The ``only if'' part}. Now suppose that a complex curve $\mathbf{X}$ is defined over $\overline{\mathbb{Q}}$. Then according to \cite{ShVo1989} it can be realized as a union of equilateral triangles, so there exists such a dessin 
$\mathcal{D}$ with $\rho_2^3\equiv1$ and the metric parameters $\ell\equiv1$,  $\varphi\equiv\frac{\pi}{3}$ that the complex structure defined by it is isomorphic to $\mathbf{X}$. Denoting $\beta_0$ the original Belyi function corresponding to the equilateral triangulation, introduce
$$
\beta:=\frac{4}{27}\frac{(\beta_0^2-\beta_0+1)^3}{\beta_0^2(1-\beta_0)^2}.
$$
It is also a Belyi function, and the transition from $\beta_0$ to $\beta$ corresponds to the {\it barycentric subdivision} of the original triangulation.

In the %corresponding
local coordinates $\{z_\varepsilon\mid\varepsilon\in\vec{E}(\mathcal{D})\}$ constructed as above: the compositional inverse to the Belyi function $\beta$ is defined by the Christoffel-Schwartz integral
$$
z_\varepsilon(P)=
\frac{\displaystyle\int\limits_0^{\beta(P)}\frac{\mathrm{d}w}{\sqrt[6]{w^5(1-w)^3}}}
{\displaystyle\int\limits_0^{1}\frac{\mathrm{d}w}{\sqrt[6]{w^5(1-w)^3}}}
$$
The ``square'' local coordinates $\{Z_\varepsilon\mid\varepsilon\in\vec{E}(\mathcal{D})\}$ are defined in terms of another
Christoffel-Schwartz integral
$$
Z_\varepsilon(P)=
\frac{\displaystyle\int\limits_0^{\beta(P)}\frac{\mathrm{d}w}{\sqrt[4]{w^3(1-w)^2}}}
{\displaystyle\int\limits_0^{1}\frac{\mathrm{d}w}{\sqrt[4]{w^3(1-w)^2}}};
$$
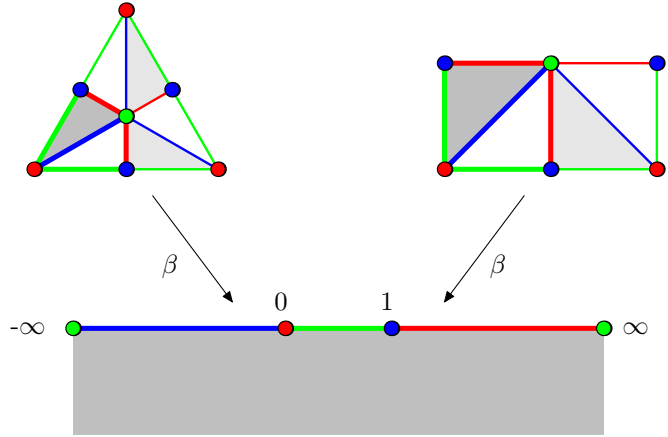
\begin{center}
\includegraphics{quad_colors-105.mps}
\end{center}

The transformation $z_\varepsilon\to Z_\varepsilon$ is defined globally on $\mathbf{X}$ and realizes the transition from the equilateral triangles to the squares. $\blacksquare$

\end{document}